\documentclass[11pt,oneside,a4paper,reqno, english]{article}

\usepackage{amsmath}
\usepackage{amssymb}
\usepackage{amsthm}
\usepackage{makeidx}
\usepackage{pifont}
\usepackage{bbding}
\usepackage{color}
\usepackage{textcomp}
\usepackage{graphicx}
\usepackage{pgf,tikz}
\date{}

\usepackage[english]{babel}
\usepackage{soul}

   \newcommand{\beq}{\begin{equation}}
   \newcommand{\eeq}{\end{equation}}
   \newcommand{\beqs}{\arraycolsep1.5pt\begin{eqnarray}}
   \newcommand{\eeqs}{\end{eqnarray}\arraycolsep5pt}
   \newcommand{\beqsn}{\arraycolsep1.5pt\begin{eqnarray*}}
   \newcommand{\eeqsn}{\end{eqnarray*}\arraycolsep5pt}

  \xdefinecolor{azzurro}{rgb}{0.2,0.9,0.9}
\xdefinecolor{violetto}{rgb}{0.5,0.5,0.9}
\xdefinecolor{marrone}{rgb}{0.8,0.4,0.2}

\def\QED{\hfill $\Box$ }
\newtheorem{pft}{}                              

\newcommand{\bpf}{\begin{pft}\begin{proof}\mbox{\bf Proof}.\ \rm}
\newcommand{\bpfof}[1]{\begin{pft}\begin{proof}{\bf Proof of #1}\
\rm}
\newcommand{\bspf}{\begin{pft}\begin{proof}{\bf Sketch of proof}.\
\rm}

\newcommand{\epf}{\nopagebreak\end{proof}\end{pft}}

\newtheorem{thm}{Theorem}[section]
\newtheorem{rem}[thm]{Remark}
\newtheorem{cor}[thm]{Corollary}
\newtheorem{prop}[thm]{Proposition}
\newtheorem{lemma}[thm]{Lemma}
\newtheorem{defn}[thm]{Definition}

\numberwithin{equation}{section}

\def \trait (#1) (#2) (#3){\vrule width #1pt height #2pt depth #3pt}
\def \qed{\hfill
        \trait (0.1) (6) (0)
        \trait (6) (0.1) (0)
        \kern-6pt
        \trait (6) (6) (-5.9)
        \trait (0.1) (6) (0)
\medskip}
\def\square{\trait (0.1) (6) (0)
        \trait (6) (0.1) (0)
        \kern-6pt
        \trait (6) (6) (-5.9)
        \trait (0.1) (6) (0)}

\def\Bbb#1{{\fam\msbfam\relax#1}}

\def\R{\Bbb R}
\def\N{{\Bbb N}}

\newcommand{\ds}{\rightarrow}

\newcommand{\la}{\lambda}

\newcommand{\om}{\Omega}
\newcommand{\Om}{\Omega}
\newcommand{\f}{\varphi}

\def\QED{\hbox{\ $\rlap{$\sqcap$}\sqcup$}\bigbreak}

\def\A{{\cal A}}
\def\B{{\cal B}}

\def\F{{\cal F}}

\def\L{{\cal L}}

\newcommand{\wi}{W^{1,\infty}(\Om)}

\def\infess{\mathop{\rm ess\: inf }}
\def\supess{\mathop{\rm ess\: sup }}

\def\dist{{\rm dist}\,}
\def\supp{{\rm supp}\,}

\newcommand{\dx}{\,dx}

\newcommand{\wto}{\rightharpoonup}
\newcommand{\weakst}{\stackrel{\ast}{\rightharpoonup}}

\newcommand{\e}{\varepsilon}

\newcommand{\li}{{L}^\infty}

\def\supess{\mathop{\rm ess\: sup }}

\font\tenmsb=msbm10 \font\sevenmsb=msbm7 \font\fivemsb=msbm5
\newfam\msbfam
\textfont\msbfam=\tenmsb \scriptfont\msbfam=\sevenmsb
\scriptscriptfont\msbfam=\fivemsb
\def\Bbb#1{{\fam\msbfam\relax#1}}

\def\rr{\Bbb R}

\newcommand {\Me}{{\cal M}}

\def\to{\rightarrow}

\def\essinf{\mbox{ess inf}}

\def\d{\, d}

\newenvironment{michelarev}{\color{blue}}{\color{black}}
\newcommand{\bmicr}{\begin{michelarev}}
\newcommand{\emicr}{\end{michelarev}}

\title{\bf  $\Gamma$-convergence  for  power-law functionals\\ with variable exponents}

\author{{\sc Michela Eleuteri}\\
Universit\`a degli Studi di Modena e Reggio Emilia\\
Dipartimento di Scienze Fisiche Informatiche e Matematiche \\
 via Campi 213/b, 41125 Modena (Italy)\\
{\tt michela.eleuteri@unimore.it}
\\[4mm] and
 {\sc Francesca Prinari}\\
 Universit\`a degli Studi di Ferrara  \\
 Dipartimento di Matematica e Informatica 
\\
 via Machiavelli 30, 44121 Ferrara (Italy)
\\
 {\tt francesca.prinari@unife.it}}

\begin{document}
\maketitle
\begin{abstract}
We study the $\Gamma$-convergence of  the functionals $F_n(u):= || f(\cdot,u(\cdot),Du(\cdot))||_{p_n(\cdot)}$  and $\F_n(u):=\displaystyle \int_{\Omega} \frac{1}{p_n(x)}  f^{p_n(x)}(x,u(x),Du(x))dx$ defined on 
$X\in \{L^1(\Omega,\R^d), L^\infty(\Omega,\R^d), C(\Omega,\R^d)\}$ (endowed with their usual norms)  
with effective domain  the Sobolev space $W^{1,p_n(\cdot)}(\Omega, \mathbb{R}^d )$. Here $\Omega\subseteq \R^N$ is a bounded open set, $N,d \ge 1$ and  the measurable functions  $p_n: \overline{\Omega} \rightarrow (1, + \infty) $   satisfy  the  conditions    $\displaystyle{\supess_{\overline \Omega} p_n  \le \, \beta \, \infess_{\overline \Omega} p_n }$  for a fixed constant $\beta > 1$ and $\displaystyle{\infess_{\overline \Omega} p_n   \rightarrow + \infty}$ as $n \rightarrow + \infty$.
We show that when  $f(x,u,\cdot)$ is  level convex  and lower semicontinuous and it satisfies a uniform growth condition  from below, then, as $n\to \infty$,   the sequences $(F_n)_n$ $\Gamma$-converges in $X$ 
 to the functional $F$ represented as  $F(u)= || f(\cdot,u(\cdot),Du(\cdot))||_{\infty}$ 
on the effective domain $W^{1,\infty}(\Omega, \mathbb{R}^d )$. Moreover we show that  the $\Gamma$-$\lim_n \F_n$   is given by the functional $ \F(u):=\left\{\begin {array}{lll} \!\!\!\!\!\!  & \displaystyle 0
&  \hbox{if  } || f(\cdot,u(\cdot),Du(\cdot)) ||_{\infty}\leq 1,\\
\!\!\!\!\!\!  & +\infty  & \hbox{otherwise in } X.\\
\end{array}\right.
$

\end{abstract}

{\small
\tableofcontents
}

\vspace{15mm}

%
%

\section{Introduction}
The  classical functionals of the Calculus of Variations 
are represented in the integral form  $$\mathcal H_n(u,A) =\int_{A} f_n(x,Du(x))dx,$$ and are defined on some subset of a Sobolev space $W^{1,q}(\Omega)$, where $A\subseteq \Omega\subseteq \R^N$ with  $A$ and $\Omega$ open sets. When the sequence of Borel functions $(f_n)$ satisfies a uniform growth condition of order $q>1,$  often named {\it standard growth condition} 
\beq\label{pcrescita}\alpha(|\xi|^q -1)\leq f_n(x,\xi)\leq \beta(|\xi|^q +1)\eeq
with $0<\alpha\leq \beta,$
 then it is possible to apply a general compactness procedure to get that there exists a subsequence $(\mathcal H_{k_n})_n $ $ \Gamma$-converging  with respect to the 
$L^q$-norm to a functional $\mathcal H_0$ that can be represented in the integral form.
\noindent If the growth condition of order $q$ is not uniformly satisfied, then the $\Gamma$-limit of the  sequence $\mathcal H_n(u,A)$  can  lose the additivity property with respect to the union of disjoint open sets and may assume a different representation form. 
This is the case, for instance, when  \eqref{pcrescita}  is replaced by the so called {\it non-standard growth condition} (considered for the first time in the pioneering papers by Marcellini \cite{M89, M91}): 
in the case of integral functionals  exhibiting a
gap between the coercivity and the growth exponent,  in \cite{MM}  Mingione and Mucci  show that  in the relaxation procedure  energy concentrations  may appear leading to a measure representation of the relaxed functional with a nonzero singular part.

A different situation appears for example in \cite{GNP}: Garroni, Nesi and Ponsiglione consider the case when \beq\label{fn}f_n(x,\xi):=\frac 1 n (f(x,\xi))^n\eeq  where   $f(x,\xi):=a(x)|\xi|$ with    $a\in L^\infty(\Omega)$ satisfying the condition $\essinf_{x\in \Omega} a(x)>0$.  This sequence   $(f_n)$ does not 
verify uniformly a $q$-growth condition and in \cite[Proposition 2.1]{GNP} it is shown that, when $\Omega$ is the unitary cube of $\R^N$,  the $\Gamma$-limit  (with respect to the $L^1$-convergence)  of the sequence $(\mathcal I_n)$ defined by 
\beq\label{Fpintro1} \mathcal I_n(u):=\left\{\begin {array}{lll} \!\!\!\!\!\! & \displaystyle   \int_{\Omega} \frac 1 n(a(x)|Du(x)|) ^ndx,
&  \hbox{if  } u\in  W^{1,n}(\Omega),\\
\!\!\!\!\!\!  &+\infty  & \hbox{otherwise in } L^1(\Omega) 
\end{array}\right.
\eeq
 is given by 
\beq\label{Ilim1}  \mathcal I(u):=\left\{\begin {array}{lll} \!\!\!\!\!\!  & \displaystyle 0
&  \hbox{if  } || a(x)|Du(x)| ||_{\infty}\leq 1,\\
\!\!\!\!\!\!  & +\infty  & \hbox{otherwise.}\\
\end{array}\right.
\eeq
Moroever, in the same paper (see Proposition 2.6 therein), it is proved  that the sequence of the $L^n$-norms  
\beq\label{Fpintro2}\ I_{n}(u):=\left\{\begin {array}{lll} 
\!\!\!\!\!\!  & \displaystyle \left(\int_{\Omega} ( a(x)|Du(x)|)^{n}dx\right)^{1/n} &  \hbox{if  } u\in  W^{1,n}(\Omega),\\
\!\!\!\!\!\!  & +\infty  & \hbox{otherwise in } L^1(\Omega). \\
\end{array}\right.
\eeq
$\Gamma$-converges (with respect to the $L^1$-convergence) to the  functional $I$  represented  in the {\sl supremal} form
\beq\label{Ilim2} I (u):=\left\{\begin {array}{lll} 
\!\!\!\!\!\!  & \displaystyle \supess_{\Omega} a(x)|Du(x)| &  \hbox{if  } u\in  W^{1,\infty}(\Omega),\\
\!\!\!\!\!\!  & +\infty  & \hbox{otherwise in } L^1(\Omega) .\\
\end{array}\right.
\eeq
Recently, the class of functionals represented in the general supremal form
 has being studied with growing interest. 
 They appear in a very natural way in variational problems  where the  relevant quantities do not express a mean property and the values of the energy densities  on  very small subsets of $\Omega$ cannot be neglected.
 Their study was introduced  by Aronsson in the
1960s (see \cite{A1}, \cite{A2}, \cite{A3}).  In the seminal papers \cite{McS} and \cite{J},
 the supremal functional $F(u)=||Du||_{\infty}$  appears in the variational problem of  finding the best Lipschitz extension $u$  in $\Omega$ of a function $g$ defined on the boundary $\partial \Omega$.  There after, several mathematical model have been formulated  by means of  a supremal functional: for example,
 models describing dielectric breakdown in a composite material (see \cite{GNP})
or polycrystal plasticity (see \cite{BN}). Also the  problem of image reconstruction and enhancement can  been formulated as an $L^\infty$ problem (see \cite{CM}). A recent application  involving a supremal functional has been given in \cite{K} where, in order to lay the rigorous mathematical foundations of the Fluorescent Optical Tomography (FOT),   the author  poses FOT as a minimisation problem in $L^\infty$ with PDE constraints. \\

\noindent The above mentioned results  contained in \cite{GNP} have been generalized in different directions:
\begin{itemize}
\item[-] \noindent when   $X=C(\bar \Omega,\R^d)$ is endowed with the uniform topology,  in \cite{CDP}, \cite{P} and \cite{PZ} the authors study the $\Gamma$-convergence of the family  of integral functional 
$F_p: X \to [0,+\infty]$ given by 
	$$ F_p(u):=\left\{\begin {array}{cl}
	\displaystyle \left( \int_{\Om} f^p(x,u(x), D u(x))dx \right)^{1/p},
	&  \hbox{if } \, u\in W^{1,p}(\Omega,\R^d),\\
	+\infty,  & \hbox{otherwise in }$X$
\end{array}\right.
$$
gradually weakening the assumptions on $f$:
\begin{itemize}
\item[-] in \cite[Theorem 3.1]{CDP} the function $f$ is a  normal integrand satisfying a superlinear  growth condition
and a  generalized Jensen inequality for gradient Young measures; in particular the $\Gamma$-convergence result  therein   holds when   the sub level sets
$\{\xi\in \R^{d\times N}:\,  f(x,u,\xi)\leq \la\}$ are closed and  convex for every $\la\in \R$;
\item[-] in \cite{P} the function $f=f(x,\xi)$ is assumed to be  a Carath\'eodory integrand satisfying the linear growth  condition \eqref{pcrescita}  with $q=1$,  while in \cite{PZ}  the continuity assumption on $f$ with respect to the gradient variable has dropped and  the function $f=f(x,\xi)$ is assumed to be only $\L^{N}\times \B_{d\times N}$-measurable;
 \end{itemize}
 \item[-] in the papers \cite{BN}  and  \cite{AP}   the space $X$ coincides with   the class of the functions $U$ in $L^1(\Omega,\R^{d\times N})$ or in  $ L^\infty(\Omega,\R^{d\times N})$, constrained to satisfy a  general rank-constant differential constraint $\mathcal{A}U=0$ and  the functionals   $F_p: X \to [0,+\infty]$ are defined  by 
	$$ F_p(U):=\left\{\begin {array}{cl}
	\displaystyle \left( \int_{\Om} f^p(x, U(x))dx \right)^{1/p},
	&  \hbox{if } \, u\in L^{p}(\Omega,\R^d)\cap X,\\
	+\infty,  & \hbox{otherwise in } X.
\end{array}\right.
$$
In \cite{BN} the authors consider the case when  $f(x,\xi)=a(x)|\xi|$ while in \cite{AP} the $\Gamma$-convergence is studied in the wider  class of   $\A_{\infty}$-quasiconvex function $f$ (see  Definition 3.2 therein).
These results have been extended in  \cite{BMP} in the setting of  variable exponent Lebesgue space when $f(x,\cdot)$ is quasiconvex in the sense of Morrey.
 \end{itemize}
A different generalization of the results contained in \cite{GNP} has been given  by  Bocea-Mihilescu  in \cite{BM}: they show   the  $\Gamma$-convergence of  the sequences  \eqref{Fpintro1} and  \eqref{Fpintro2} respectively to  the functionals $\mathcal I$ and $ I$  given by \eqref{Ilim1} and \eqref{Ilim2} respectively
when $\Omega$ is a Lipschitz connected open set satisfying $\mathcal{L}^N(\Omega)=1$ and  the constant sequence $(n)$ is replaced by a sequence   $(p_n)=(p_n(x))$ of  Lipschitz continuous functions satisfying
$$p_n^- :=  \infess_{\overline \Omega} p_n   \rightarrow + \infty
$$ as $n \rightarrow + \infty$ and 
$$ p_n^+=\supess_{\overline \Omega} p_n  \le \, \beta  p_n^- $$
  for a fixed constant $\beta > 1$.

Inspired by  \cite{BM}, in our paper we consider the general case   when $f:\Omega\times \R^d\times \R^{N\times d}\to [0,+\infty)$    in  formula \eqref{fn} is a Borel function  satisfying the  coercivity assumption  \beq\label{basso} f(x,u, \xi) \geq  \alpha |\xi|^{\gamma}  
\qquad \hbox{ for a.e } x\in \Om, \hbox{ for every }
(u, \xi)\in\R^d\times  \R^{Nd}
\eeq
 (with $\alpha,\gamma>0$) 
 and such that  its sub level sets $\big\{\xi\in\R^N\colon f(x,u,\xi)\le t\big\}$ are closed and convex for any $t\in \R$.
 The last assumption is  sufficient  to ensure the lower semicontinuity  with respect to the weak* lower semicontinuity of the supremal functional 
$$ F(u):= \supess_{\Omega} f(x,u(x),Du(x))$$
in the space $W^{1,\infty}(\Omega,\R^d)$ (see Theorem 3.4 in \cite{BJW}).

Under the previous hypotheses,  in Theorems  \ref{gamma} and \ref{gamma2} we show that,   if we consider $X\in \{L^1(\Omega,\R^d), L^\infty(\Omega,\R^d), C(\Omega,\R^d)\},$
 endowed  with their  usual norms,  
then  the sequence of  functionals  $F_n:X\to [0,+\infty]$ 
defined by 
$$
F_n(u):=\left\{\begin {array}{cl} \displaystyle \ || f(\cdot,u(\cdot),Du(\cdot))||_{p_n(\cdot)}
&  \hbox{if } \, u\in W^{1,p_n(\cdot)}(\Omega,\R^{d})\\
+\infty  & \hbox{otherwise,}
\end{array}\right.
$$
as $n\to +\infty$ $\Gamma$-converges to the functional $F$ given  by  
$$F(u):=\left\{\begin {array}{cl} \displaystyle
\supess_{\Omega} f(x,u(x),Du(x))
&  \hbox{if } \, u\in W^{1,\infty}(\Omega,\R^{d}),\\
+\infty  & \hbox{otherwise in } X.
\end{array}\right.
$$
In particular, thanks to the more general growth condition \eqref{basso} on $f$,  we get an improvement of Theorem 3.1 in \cite{CDP} (see Corollary \ref{gamma4}).
 
Moreover in Theorems \ref{gamma5} and \ref{gamma6} we show that the sequence of the integral functionals  $\F_n:X\to [0,+\infty]$  defined by 
$$
\F_n(u):=\left\{\begin {array}{cl} \displaystyle \ \int_{\Omega} \frac{1}{p_n(x)}  f^{p_n(x)}(x,u(x),Du(x))dx
&  \hbox{if } \, u\in W^{1,p_n(\cdot)}(\Omega,\R^{d})\\
[2mm]
+\infty  & \hbox{otherwise,}
\end{array}\right.
$$
as $n\to +\infty$, $\Gamma\hbox{-}$converges
 to the functional  defined by 
$$ \F(u):=\left\{\begin {array}{cl} \displaystyle 0
&  \hbox{if } \, u\in W^{1,\infty}(\Omega,\R^{d}) \hbox{ and }||f(x,u(x),Du(x))||_{\infty}\leq 1,\\
+\infty  & \hbox{otherwise in } X.
\end{array}\right.
$$
  Note that  the proofs of  the previous results  are given in the vectorial setting and in general case $ \L^N(\Omega)\in (0,+\infty)$.
  Moreover we do not need any connectedness hypothesis  on $\Omega$ and only when  $X=L^1(\Omega,\R^d)$ we assume that $\partial\Omega$ is Lipschitz regular.
We  point out that in \cite{BM}  the  weak lower semicontinuity  in $L^q(\Omega)$  of  the integral functionals \eqref{Fpintro1} and \eqref{Fpintro2}    was enough to prove the $\Gamma$-convergence results therein. 
Instead, the more general class of our variational functionals requires  as key tool the use of {\it gradient Young measures}: indeed, this instrument turns to be crucial in order to show the $\Gamma$-liminf inequality,  combined with  a Jensen type inequality satisfied  by level convex functions, see Theorem \ref{Jensen}.  Moreover, in our paper, we  deal with more general topologies instead of treating only with the strong convergence in $L^1(\Omega)$ as in \cite{BM}. This allows us to remove the regularity assumptions on  $\partial \Omega$ in the case of $X = \{L^\infty(\Omega,\R^d), C(\Omega,\R^d)\}.$\\
\noindent We devote a forthcoming paper  to study  the homogenization of supremal functionals of the form
$$F_{\epsilon}(u):= \supess_{\Omega}  g\left(\frac x \epsilon, Du(x)\right)$$
where $g(x,\xi):=(f(x, \xi))^{p(x)}.$
With this  aim, we   will proceed our analysis by  discussing the $\Gamma$-convergence of the sequence of functionals  
$$
H_n(u):=\left\{\begin {array}{cl} \displaystyle \left( \int_{\Omega} \frac{1}{n p(x)}  f^{np(x)}(\cdot,u(\cdot),Du(\cdot))dx\right)^{1/n}
&  \hbox{if } \, u\in W^{1,p_n(\cdot)}(\Omega,\R^{d})\\
[2mm]
+\infty  & \hbox{otherwise,}
\end{array}\right.
$$
already  considered  in the case  $f(\xi)=|\xi|$ by Zhikov in \cite{zh1} (see also \cite{zh2}) and, more recently, by Bocea and Mihilescu in \cite{BM} in the case  when  $f(x,\xi)=a(x)|\xi|$.

\section{Preliminary results}

\subsection{Variable exponents Lebesgue-Sobolev spaces}
In  this section we collect some basic results concerning variable exponent Lebesgue and Sobolev spaces. For more details we refer to the monograph \cite{DHHR11}, see also \cite{KR91}, \cite{ELN99}, \cite{ER00}, \cite{ER02}.
\\
For the purpose of our paper, we consider the case when  $\Omega \subset \mathbb{R}^N$ is  an open set (where $N\geq 1$) and denote by $\mathcal{L}^N(\Omega)$ the $N$-dimensional Lebesgue measure of the set $\Omega$.  In the sequel we will consider functions $u: \Omega \rightarrow \mathbb{R}^d$, with $d \ge 1$ and we denote by $k$ any dimension different from $Nd$. 
\\
For any Lebesgue measurable function $p: \Omega \rightarrow [1, + \infty]$
we define
$$
p^- := \displaystyle \infess_{x \in \Omega} p(x) \qquad \qquad p^+ :=\displaystyle
\supess_{x \in \Omega} p(x).
$$
Such function $p$ is called {\it variable exponent} on $\Omega$. If $ p^+<+\infty$  then we call $p$ a {\it bounded variable exponent}.
\\
In the sequel we need to introduce the Lebesgue spaces with variable exponents, $L^{p(\cdot)}$. They differ from the classical $L^p$ spaces because now the exponent $p$ is not constant but it is a variable exponent in the sense specified above. Originally the spaces $L^{p(\cdot)}$  have been introduced in the case $1 \le p^- \le p^+ < \infty$ by Orlicz \cite{O31} in 1931 and,  in the case $p^+ = \infty$,   by Sharpudinov \cite{S79} and later  (in the higher dimensional case), by Kov\'a\v{c}ik and R\'akosn\'ik \cite{KR91}.

In the sequel we consider the case $p^+ < \infty$. In this case the variable exponent Lebesgue space $L^{p(\cdot)}(\Omega)$ can be defined as \[
L^{p(\cdot)}(\Omega) := \left \{u: \Omega \rightarrow \R \,\,\, \textnormal{measurable such that} \,\, \int_{\Omega} |u(x)|^{p(x)} \, dx < + \infty   \right \}.
\]

Let us note that in the case  $p^+ = \infty$ the  space above defined  may even fail to be a vector space (see \cite{CM02} Section 2)
and a different definition of the variable Lebesgue spaces  has been  given in order to preserve the vectorial structure of the space (we refer to \cite{DHHR11}, Definition 3.2.1).
On the other hand if $p^+ < + \infty,$ it is possible to show that $L^{p(\cdot)}(\Omega)$ is a Banach space endowed with the {\it Luxemburg norm} 
\[
\|u\|_{p(\cdot)} := \inf \left \{\lambda  > 0: \,\,\, \int_{\Omega} \left |\frac{u(x)}{\lambda} \right |^{p(x)} \, dx \le \, 1  \right \}
\]
(see  Theorem 3.2.7 in  \cite{DHHR11}). 
Moreover if $p^+ < + \infty$ the space $L^{p(\cdot)}(\Omega)$ is separable and the space  $\mathcal{C}^{\infty}_0(\Omega)$ is dense in $L^{p(\cdot)}(\Omega)$  while if $1 < p^- \le p^+ < + \infty$ the space $L^{p(\cdot)}(\Omega)$ is reflexive and uniformly convex (see 
Lemma 3.4.1, Theorem 3.4.7, Theorem 3.4.9 and Theorem 3.4.12  in \cite{DHHR11}). Finally, by 
Corollary 3.3.4 in \cite{DHHR11}, if $0 <  \mathcal{L}^N(\Omega) < + \infty$ and $p$ and $q$ are variable exponents such that $p \le q$   a.e. in $\Omega$, then the embedding $L^{q(\cdot)}(\Omega) \hookrightarrow L^{p(\cdot)}(\Omega)$ is continuous.  The embedding constant is less or equal to  
$2  \max  \left \{  {\mathcal{L}^N(\Omega)}^{{(\frac 1 q-\frac 1 p)}^+},  {\mathcal{L}^N(\Omega)}^{{(\frac 1 q-\frac 1 p)}^-}  \right \}.$

For any variable exponent $p$, we define $p'$ by setting
\[
\frac{1}{p(x)} + \frac{1}{p'(x)} = 1,
\]
with the convention that, if $p(x) = \infty$ then $p'(x) = 1$. The function $p$ is called {\it the dual variable exponent of $p$}.

We have the following result (for more details see Lemma 3.2.20 in \cite{DHHR11}).

\begin{thm} {\sl (H\"older's inequality)}
Let $p,q,s$ be measurable exponents such that
\[
\frac{1}{s(x)} = \frac{1}{p(x)} + \frac{1}{q(x)}
\]
a.e. in $\Omega$. Then
\[
\|fg\|_{s(\cdot)} \le \, \left( \left(\frac{s}p\right)^+ + \left(\frac{s}q\right)^+ \right)\, \|f\|_{p(\cdot)} \, \|g\|_{q(\cdot)}
\]
for all $f \in L^{p(\cdot)}(\Omega)$ and $g \in L^{q(\cdot)}(\Omega)$, where in the case $s = p = q = \infty$, we use the convention $\frac{s}{p} = \frac{s}{q} = 1.$
\\
In particular, in the case $s  = 1$, we have
\[
\left | \int_{\Omega} f \, g \, dx\right | \le \, \int_{\Omega} |f| \, |g| \, dx \le {\, \left( \frac{1}{p^-} + \frac{1}{p'^-} \right)}\, \, \|f\|_{p(\cdot)} \, \|g\|_{p'(\cdot)}.
\] 
\end{thm}

We moreover introduce the {\it modular} of the space $L^{p(\cdot)}(\Omega)$ which is the mapping $\rho_{p(\cdot)}: L^{p(\cdot)}(\Omega) \rightarrow \mathbb{R}$ defined by
\[
\rho_{p(\cdot)}(u) := \int_{\Omega} |u(x)|^{p(x)} \,dx.
\]
Thanks to Lemma 3.2.4  in \cite{DHHR11}, for every $u \in L^{p(\cdot)}(\Omega)$ 
\begin{eqnarray}
&& \|u\|_{p(\cdot)}  \leq 1 \Longleftrightarrow \rho_{p(\cdot)}(u)\leq  1 \label{mod-1}\\
&& \|u\|_{p(\cdot)}  \leq 1 \Longrightarrow \rho_{p(\cdot)}(u)\leq    \|u\|_{p(\cdot)}\label{mod-2}\\
&& \|u\|_{p(\cdot)}  > 1 \Longrightarrow \|u\|_{p(\cdot)}\le \, \rho_{p(\cdot)}(u)  \label{mod-3}.
\end{eqnarray}

The following further results hold in the special case $p^+<+\infty$.
By Lemma 3.2.5 in \cite{DHHR11},    for every $u \in L^{p(\cdot)}(\Omega)$ it holds
\begin{equation} \label{relaztotale}\min \left\{\big( \rho_{p(\cdot)}(u)\big)^{\frac 1 {p^-} } , \big(\rho_{p(\cdot)}(u)\big)^{\frac 1 {p^+} }  \right \} \leq  \|u\|_{p(\cdot)} \leq \max  \left \{  \big( \rho_{p(\cdot)}(u)\big)^{\frac 1 {p^-} } , \big(\rho_{p(\cdot)}(u)\big)^{\frac 1 {p^+} } \right \}.
\end{equation}
In particular, we get that
\beq\|1\|_{p(\cdot)} \leq \max\{\big(\L^N(\Omega) \big)^{\frac 1 {p^-} } , \big(\L^N(\Omega)\big)^{\frac 1 {p^+} }   \} . \label{mod4}\eeq
Moroever, from \eqref{relaztotale}, taking into account  \eqref{mod-2}, \eqref{mod-3}, it follows that for every $u \in L^{p(\cdot)}(\Omega)$  \begin{eqnarray}
&& \|u\|_{p(\cdot)}  > 1 \Longrightarrow \|u\|_{p(\cdot)}^{p^-} \le \, \rho_{p(\cdot)}(u) \le \, \|u\|^{p^+}_{p(\cdot)} \label{mod1}\\
&& \|u\|_{p(\cdot)} < 1 \Longrightarrow \|u\|_{p(\cdot)}^{p^+} \le \, \rho_{p(\cdot)}(u) \le \, \|u\|^{p^-}_{p(\cdot)} \label{mod2}.
\end{eqnarray}
Finally, thanks to Lemma 3.4.2  in \cite{DHHR11},   for every $u \in L^{p(\cdot)}(\Omega)$
\begin{eqnarray}
&& \|u\|_{p(\cdot)}  < 1 \Longleftrightarrow \rho_{p(\cdot)}(u)<  1 \label{mod0}\\
&&\|u\|_{p(\cdot)}  = 1 \Longleftrightarrow \rho_{p(\cdot)}(u)  = 1 \label{mod3}\\
&&\|u\|_{p(\cdot)}  > 1 \Longleftrightarrow \rho_{p(\cdot)}(u)  >1 \label{mod5}.
\end{eqnarray}

%
%
%

 We conclude this part by recalling the definition of  variable exponent Sobolev spaces.  For more details we refer to \cite{CM02} (see also \cite{DHHR11}, Definition 8.1.2). 

\begin{defn}
 
Let $k,d\in \N$, $k\geq 0$,  and let  $p$ be a measurable exponent. We define
$$W^{k, p(\cdot)}(\Omega,\R^d):=\{ u:\Omega\to \R^d :u, \partial_{\alpha} u \in L^{p(\cdot)}(\Omega,\R^d) \quad \forall \alpha \hbox{  multi-index such that $|\alpha| \le \, k$ }  \},$$
where $$L^{p(\cdot)}(\Omega,\R^d):=\{ u:\Omega\to \R^d :\ |u| \in L^{p(\cdot)}(\Omega)  \}.$$
We define the semimodular on $W^{k, p(\cdot)}(\Omega)$ by
\[
\rho_{W^{k, p(\cdot)}(\Omega)}(u) := \sum_{0 \le |\alpha| \le k} \rho_{L^{p(\cdot)}(\Omega)}(|\partial_{\alpha} u|)
\]
which induces a norm  by 
\[
\|u\|_{W^{k, p(\cdot)}(\Omega)} := \inf \left \{ \lambda > 0: \,\, \rho_{W^{k, p(\cdot)}(\Omega)}\left (\frac{u}{\lambda} \right ) \le \, 1 \right \}.
\]
\end{defn}
 For $k \in \mathbb{N} \setminus \{0\},$ the space $W^{k, p(\cdot)}(\Omega)$ is called  {\it Sobolev space} and its elements are called {\it Sobolev functions}. Clearly $W^{0, p(\cdot)}(\Omega) = L^{p(\cdot)}(\Omega)$.

\subsection{Level convex functions}

\begin{defn}We say that $f:\R^k\to \R$ is level convex if
 for every $t\in\R$ the level set $\big\{\xi\in\R^k\colon f(\xi)\le t\big\}$ is convex. 
\end{defn}
We recall Jensen's inequality introduced by  Barron, Jensen, and Liu in \cite{BJL}  for lower semicontinuous and level convex functions (see also \cite{BJW} Theorem 1.2).  \begin{thm}\label{Jensen}  
Let $f:\R^k \to \R$ be a lower semicontinuous and level convex
function, and let $\mu$ be a probability measure supported on the open set  $\Omega\subseteq \R^N$. Then for every function  $u\in L^1_{\mu}(\Omega;\R^k)$ we have
\begin{equation} \label{JI} f\left(\int_{\Omega} u(\xi)  d\mu(\xi)\right)\le \mu\hbox{-}\supess_{\xi\in \Omega} (f\circ u)(\xi).
\end{equation}
\end{thm}

\subsection{Young measures}\label{YMe}
In this section we briefly recall some results on the theory of Young measures (see e.g. \cite{Ba89}, \cite{BL}). If $\Omega\subseteq \R^N$ is an open set (not necessarily
bounded) and  $d\geq 1$, we denote by $C_c(\Omega;\R^d)$ the set of continuous
functions with compact support in $\Omega$, endowed with the supremum
norm. The dual of the closure of $C_c(\Omega;\R^d)$  may be identified
with the set of $\rr^d$-valued Radon measures with finite mass
$\Me(\Omega;\rr^d)$, through the duality
$$
\langle \mu, \varphi\rangle := \int_\Omega \varphi(\xi)\d\mu (\xi)\,, \qquad
\mu\in \Me(\Omega;\rr^d)\,,\qquad \varphi\in C_c(\Omega;\rr^d)\,.
$$

\begin{defn} A map $\mu:\Omega\mapsto \Me(\Omega;\rr^d)$ is said to be
weak$^*$-measurable if $x\mapsto \langle \mu_x, \varphi\rangle$ are
measurable for all $\varphi\in C_c(\Omega;\rr^d)$.
\end{defn}

\begin{defn} Let $(V_n)$  be a bounded sequence in  $L^1(\Omega,\R^d)$.
We say that $(V_n)$ is {\it equi-integrable} if for all $\e>0$
there exists $\delta>0$ such that, for every measurable  $E\subset \Omega,$ if $\L^N(E)<\delta,$ then
$$
\sup_n\int_E |V_n(x)|\dx<\e\,.
$$
For every $1<p<+\infty$ we say that $(V_n)$ is $p$-{\it equi-integrable} if $(|V_n|^p)$
is equi-integrable.
\end{defn}

We are in position to state the main result concerning Young Measures (for a proof see \cite[Theorem 3.1]{Mu}).

\begin{thm}[\sl Fundamental Theorem on Young Measures]\label{FTYM}
Let $E  \subseteq  \rr^N$ be a measurable set of finite measure and let
$(V_n)$ be a sequence of measurable functions, $V_n: E\mapsto
\R^d$. Then there exists a subsequence $(V_{n_k})$ and a weak$^*$-measurable map $\mu:E\mapsto \Me(\Omega; \R^d)$ such that the following  statements hold: 
\begin{enumerate}
\item[(1)] $\mu_x\ge 0$,
$\displaystyle{\Vert\mu_x\Vert_{\Me(\Omega; \R^d)}=\int_{\R^d}d\mu_x\le 1}$ for a.e. $x\in
E$;
\item[(2)]  $\forall \f\in C_c(\Omega; \R^d)$ 
$$\f(V_{n_k})\weakst  \bar \f$$ where 
$$
\bar{\f}(x):= \langle \mu_x , \f\rangle; \qquad  {\rm for\: a.e.}\: x\in E\,;$$
\item[(3)]  for every compact subset  $K\subset \R^d$,  if $\dist (V_{n_k}
K)\to 0$ in measure, then
$$
\supp \mu_x\subset K \qquad  {\rm for\: a.e.}\: x\in E\,;
$$
\item[(4)]  $\Vert\mu_x\Vert_{\Me(\R^d)}=1$ for a.e. $x\in E$ if
and only if
\beq\label{tight}
\lim_{M\to \infty} \sup_k \L^N(\{|V_{n_k}|\ge M\})=0\,;
\eeq
\item[(5)]  if $\Vert\mu_x\Vert_{\Me(\R^d)}=1$ for a.e. $x\in E$  then in   (3) we may replace ''if'' with ``if and only if'';
\item [(6)]  if $\Vert\mu_x\Vert_{\Me(\R^d)}=1$ for a.e. $x\in E$ and $A\subseteq E$ is measurable and if  $ \f\in C(\R^d)$ is such that $ (\f(V_{n_k}))$ is  equi-integrable in  $L^1(A,\R^d)$ then
$$\f(V_{n_k}) \wto  \bar \f.$$ 
\end{enumerate}
\end{thm}

The map $\mu:E\mapsto \Me(\R^d)$ is called {\it Young measure generated by the sequence} $(V_{n_k})$. \\


%
%

From now on, for the sake of simplicity, we denote by $(V_n)$ the sequence $(V_{n_k})$ generating the corresponding Young measure. 
 
\begin{rem}\label{rem2.3} {\rm
Condition \eqref{tight} holds if there exists any $q\geq 1$ such that
$$\sup_{n\in \N} {\|V_n\|_q}<+\infty.$$
Indeed, by Chebyshev's inequality, 
$$\sup_{n\in \N} \L^N(\{|V_{n}|\ge M\})\leq \frac 1 {M^{q}} \int_{E} |V_{n}|^q dx\leq \frac C {M^q}.$$

 In particular, if  $(V_n)$   an equi-integrable sequence in  $L^1(E,\R^d)$, as a consequence of Theorem \ref{FTYM}(4)-(6) (with $\f=Id$), it generates the Young measure $\mu=(\mu_x)$ satisfying   $\Vert\mu_x\Vert_{\Me(\R^d)}=1$  such that 
$ V_{n_k} \wto  \bar{V} \hbox{ weakly in } L^1(E,\R^d)\,
$, where $$
\bar{V}(x)=  \int_{\R^d} \xi \d
\mu_x(\xi) \qquad  {\rm for\: a.e.}\: x\in E\,.
$$

%
%
}\end{rem}
The following Corollary \ref{2var} allows us to treat limits of  integrals in the form $\int_E f(x, V_{n}(x), DV_{n}(x))dx$ without any convexity assumption of $f(x,u,\cdot)$ (see Corollary 3.3 in \cite{Mu}).
\noindent First we recall the following definitions. 
\begin{defn} A function $f:\Om\times \R^{k} \to\R$ is said to be a normal integrand if  
\begin{itemize}
\item [-] $f$ is  $\L^N\otimes \B_k$-measurable;
\item  [-]  $f(x,\cdot)$ is lower semicontinuous  for a.e. $x\in \Omega$;
\end{itemize}
A function  $f:\Om\times\R^d\times \R^{k}\to\R$ is said to be a normal integrand if  
\begin{itemize}
\item  [-]  $f$ is  $\L^N\otimes\B_d\otimes \B_{k}$-measurable;
\item  [-] $f(x,\cdot,\cdot)$ is lower semicontinuous  for a.e. $x\in \Omega$;
\end{itemize}
\end{defn}
\begin{cor} \label{2var}Suppose that the sequence of measurable functions $V_{n}: E\mapsto
\R^d$ generates the Young measure $(\mu_x)$. 
\begin{enumerate} 
\item If $f:E\times \R^d\mapsto\rr$ is a  normal integrand such that  the negative part $f(x, V_{n}(x))^-$ is weakly relatively compact in $L^1(E,\R^d),$ then 
$$
\liminf_{n\to \infty} \int_E f(x, V_{n}(x))\dx \ge
\int_E \bar{f}(x)\dx\,,
$$

where
$$
\bar{f}(x):= \langle \mu_x , f(x,\cdot)\rangle = \int_{\R^d} f(x,y)\d
\mu_x(y)\,.
$$
\item if $f$ is a Carath\'eodory integrand such that  $(|f(\cdot, V_{n}(\cdot))|)$ is equi-integrable, then
$$
\lim_{n\to \infty} \int_E f(x, V_{n}(x))\dx = \int_E
\bar{f}(x)\dx<+\infty.
$$
\end{enumerate}
\end{cor}

\begin{rem}\label{limitedebole} {\rm  If   $\Omega\subseteq \R^N$ is a bounded open set and  $u_n\wto u$ in $W^{1,p}(\Omega,\R^d)$,   then the sequence  $(Du_n)_n$ is equi-integrable and generates a Young measure  $\mu=(\mu_x)$ such that $\Vert\mu_x\Vert_{\Me(\R^d)}=1$ for a.e. $x\in \Omega$  and 
$$ Du(x)= \int_{\R^{Nd}} \xi  \d
\mu_x(\xi)\,.
$$

Such a Young measure $\mu$ is usually called a $W^{1,p}$-{\sl gradient Young measure}, see \cite{Pe}.\\
Moreover, by Corollary 3.4 in \cite{Mu},  the couple $(u_n,Du_n)$ generates the Young measure $x\to \delta_{u(x)}\otimes \mu(x)$, and, if $f:\Om\times\R^d\times\R^{Nd}\to\R$ is a normal integrand bounded from below then, by Corollary \ref{2var} (1),  it follows that
\beq\label{limsotto}
\liminf_{n\to \infty} \int_\Omega f(x, u_n(x),Du_{n}(x))\dx \geq  \int_\Omega\int_{\R^{Nd}}
f(x,u(x),\xi)d\mu_{x}(\xi)   \dx.
\eeq
}\end{rem}
\subsection{$\Gamma$-convergence}
We recall the sequential characterization of  the $\Gamma$-limit when $X$ is a metric space.

\begin{prop}[\cite{DM93} Proposition 8.1]\label{seqcharac}
Let $X$ be a metric  space and let $\f_n: X \ds \R \cup
\{\pm \infty\}$ for every $n\in \N$.
 Then  $(\f_n)$ $\Gamma$-converges to $\f$ with respect to the strong topology of  $X$ (and we write $ \Gamma(X)\hbox{-}\lim_{n\to \infty}\f_n=\f$)  if and only if
\begin{description}
\item [(i)] {\rm($\Gamma$-liminf inequality)} for every $x\in X$ and for every sequence $(x_n)$ converging
to $x$, it is
$$  \f(x)\le \liminf_{n\to \infty} \f_n(x_n);$$
\item [(ii)]{\rm ($\Gamma$-limsup inequality)} for every $x\in X,$ there exists a sequence $(x_n)$  converging
to $x \in X$ such that
$$  \f(x)=\lim_{n\to \infty} \f_n(x_n).$$
\end{description}
\end{prop}
We recall that the  $\Gamma\hbox{-}\lim_{n\to \infty}\f_n$ is lower semicontinuous on $X$ (see \cite{DM93} Proposition 6.8).

Finally we recall also  that the function  $\f=\Gamma(w^*\hbox{-}X)\hbox{-}\lim_{n\to \infty}\f_n$ is weakly* lower semicontinuous  on $X$ (see \cite{DM93} Proposition 6.8) and when  $\f_n=\psi$  $\forall n\in\N$ then   $\f$   coincides with the  weakly*  lower semicontinuous (l.s.c.) envelope  of $\psi$, i.e. 
\beq\label{rilassato}
\f(x)=\sup\big\{h(x):\  \forall\, h:X\to\R\cup\{\pm\infty\} \   \ w^* \hbox
{ l.s.c.}, \  h\le \psi\hbox { on } X\big\}
\eeq
(see Remark 4.5 in \cite{DM93}).
\\
We will say that a family $(\f_p)$ $\Gamma$-converges to $\f$, with respect to the topology considered on $X$ as $p\to \infty$,
if $(\f_{p_n})$ $\Gamma$-converges to $\f$ for all sequences $(p_n)$ of positive numbers diverging to $\infty$ as $n\to\infty$.

\smallskip

\noindent Finally  we state the fundamental theorem of $\Gamma$-convergence.
\begin{thm}\label{convmin} 
Let $(\f_n)$ be an equi-coercive sequence $\Gamma$-converging on $X$ to the function $\f$ with respect to the topology of  $X$. Then we have the convergence of minima
$$\min_X \f=\lim_{n \to \infty} \inf_X \f_n.$$
Moreover we have also the convergence of minimizers: if $(x_n)$ is such that  $\lim_{n\to\infty}\f_n(x_n)= \lim_{n\to\infty}\inf_X \f_n$ then, up to subsequences, $(x_n)\to x$ and $x$ is a minimizer for $\f$.
\end{thm}
\noindent For an introduction to $\Gamma$-convergence we refer to the books \cite{DM93} and \cite{B}.

\section{Some technical lemmas}

We devote this section to show some auxiliary results  necessary in order to prove the main theorems of this paper. First of all we recall the following lemma (see  \cite[Lemma 3.2.6]{DHHR11}). We need to assume $s<p^-$ in order to ensure that $\frac {p(x)} s \geq 1 \hbox{ for a.e. } x\in \Omega.$ 
\begin{lemma}
\label{lemma3bocea}
Let $\Omega \subset \mathbb{R}^N$ be an open set  and let $p: \Omega \rightarrow [1, + \infty)$ be a bounded variable exponent. Then
\[
\||u|^s\|_{\frac{p(\cdot)}{s}}^{1/s} = \|u\|_{p(\cdot)} \qquad \textnormal{for all $u \in L^{p(\cdot)}   (\Omega)$ and $s \in (1, p^-)$.}
\]
\end{lemma}
\proof By definition
\begin{eqnarray*} 
\||u|^s\|_{\frac{p(\cdot)}{s}} &=& \inf \left \{\lambda > 0: \,\,\, \int_{\Omega} \left( \frac{|u(x)|^s}{\lambda} \right)^{\frac {p(x)}{s}} \, dx \le \, 1  \right \}\\
&=&\inf \left \{\la^s> 0: \,\,\, \int_{\Omega} \left( \frac{|u(x)|}{\la} \right)^{ {p(x)}} \, dx \le \, 1  \right \}= \|u\|_{p(\cdot)}^s.
\end{eqnarray*} 
\qed
\begin{lemma}\label{pmenoq} Let $p: \Omega \rightarrow [1, + \infty)$ be a bounded variable exponent such that  $\L^N(\Omega)<+\infty.$  Assume that there exists $\beta > 1$ such that $p^+ \le \beta p^-$.  Then for every $1\leq q\leq p^-$  and for every  $u\in  L^{p(\cdot)}(\Omega)$.
$$ ||u||_{q}\leq \max\bigg\{\big(\L^N(\Omega )\big)^{\frac {1} {q} -\frac {1} {p^-} } , \big(\L^N(\Omega)\big)^{\beta \big (\frac {1} {q}-\frac {1} {p^+}\big) }  \bigg\}  \bigg[ 1+\frac q {p^{+}(\cdot)} (\beta- 1)\bigg ]^{1/q}||u||_ {p(\cdot)}. $$ 
In particular, if $u\in L^{p(\cdot)}(\Omega,\R^d),$ then $u\in L^q(\Omega,\R^d)$ for every $1\leq q\leq p^-$. 
\end{lemma}

{\bf Proof.}
 Let $q\geq 1$. By using H\"older's inequality we have that 
$$ \int_{\Omega }|u(x)|^q dx\leq \bigg[ \frac 1 {(\frac {p(\cdot)}{ p(\cdot)- q})^{-}} + \frac 1 {(\frac {p(\cdot)}{ q})^{-}}    \bigg] ||1||_{\frac {p(\cdot)}{ p(\cdot)- q}} |||u|^q||_{\frac {p(\cdot)} {q} }=\bigg[ \frac { p^{-} - q}{ p^{+}}+  \frac{ q}    {p^{-}} \bigg]||1||_{\frac {p(\cdot)}{ p(\cdot)- q}} |||u|^q||_{\frac {p(\cdot)} {q} }. $$
 By \eqref{mod4} 
$$\|1\|_{_{\frac {p(\cdot)}{ p(\cdot)- q}}} \leq \max\left\{\big(\L^N(\Omega )\big)^{\frac {p^--q} {p^+} } , \big(\L^N(\Omega)\big)^{\frac {p^+ -q} {p^-} }  \right\}\leq \max \left \{\big (\L^N(\Omega)\big)^{1-\frac {q} {p^-} } , \big(\L^N(\Omega)\big)^{\beta \big (1-\frac {q} {p^+}\big) } \right  \}=:C $$
and 
 we get that
\begin{eqnarray*} \int_{\Omega }|u(x)|^q dx&\leq& C \bigg( \frac { p^{+} - q}{ p^{+}}+  \frac{ q}    {p^{-}} \bigg) || |u|^q ||_{\frac {p(\cdot)} {q} }\\
&=&
C\bigg[ 1+\frac q {p^{+}} \bigg (\frac  {p^{+}}{ p^{-}}- 1\bigg)\bigg ] |||u|^q ||_{\frac {p(\cdot)} {q} }\\
 &\leq &C \bigg[ 1+\frac {q} {p^{+}} (\beta- 1) \bigg ] |||u|^q ||_{\frac {p(\cdot)} {q}},
  \end{eqnarray*}
where we used the fact that $p^+ \le \, \beta \, p^-;$ this in turn implies 
\begin{eqnarray*}  ||u||_{q}\leq C^{1/q} \bigg[ 1+\frac q {p^{+}} (\beta- 1)\bigg ]^{1/q}||  |u|^q ||_{\frac {p(\cdot)} {q} }^{1/q}. \end{eqnarray*}
By Lemma \ref{lemma3bocea}, we get 
$$ ||u||_{q}\leq \max\left \{\big(\L^N(\Omega )\big)^{\frac {1} {q} -\frac {1} {p^-} } , \big(\L^N(\Omega)\big)^{\beta \big (\frac {1} {q}-\frac {1} {p^+}\big) }  \right \}  \bigg[ 1+\frac q {p^{+}} (\beta- 1)\bigg ]^{1/q}||u||_ {p(\cdot)}. $$ \qed

 In  \cite[Lemma 2]{BM}, by assuming  that $\L^N(\Omega)=1$ the sequence $\{p_n\}$ of functions $p_n: \overline{\Omega} \rightarrow (1, + \infty) $ satisfies  the  conditions:
\begin{eqnarray}
&& p_n^- \rightarrow + \infty \qquad \textnormal{as $n \rightarrow + \infty$} \label{pn1}\\[2mm]
&& \exists \, \beta > 1 : \,\,\, p_n^+ \le \, \beta \, p_n^-  \quad   \forall n \in \mathbb{N},\label{pn2} 
\end{eqnarray}
  the authors   show that if $u\in L^\infty(\Omega),$ then the $L^\infty$-norm is the limit of the $L^{p_n(\cdot)}$-norms.
 We improve   their result by showing that if the limit of the $L^{p_n(\cdot)}$-norms is finite, then $u\in L^\infty(\Omega)$. For sake of completeness, we give the detailed proof   when $\L^N(\Omega)\in (0,+\infty)$.

\begin{prop}\label{infty}
Assume $\L^N(\Omega)<+\infty$ and  let  $u:\Omega\to  \bar \R$ be a measurable function. If $(p_n)$ satisfies \eqref{pn1} and \eqref{pn2}, then  the following properties are equivalent:
\begin{description}
\item [(i)] $u\in L^\infty(\Omega)$;
\item [(ii)] $\lim_{n\to \infty} ||u||_{p_n(\cdot)}\in \R.$
\end{description}
Moreover if (i) or  (ii) holds,  then
 \beq\label{normainfty}||u||_{\infty}=\lim_{n\to \infty} ||u||_{p_n(\cdot)}.\eeq
\end{prop}
{\bf Proof.}
 \fbox{$(i) \Longrightarrow (ii) $}
  Note that,  in order to show (\ref{normainfty}), it sufficient to prove that 
$$\lim_{n\to \infty} ||u||_{p_n(\cdot)} =1 \quad \forall u\in L^\infty(\Omega) \hbox{ s.t. } ||u||_{\infty}=1.$$
\\
Let $u\in L^\infty(\Omega)$ such that  $||u||_{\infty}=1$. Then for every $n\in \N$ we  get   
 \beq\label{interm3}
\int_{\Omega }|u(x)| ^{p_n(x)}dx\leq \, \L^N(\Omega).
\eeq Since $|u(x)|\leq 1$ for a.e. $x\in \Omega$,  we have that  for every $n\in \N$
$$   \int_{\Omega} |u(x)| ^{p_n(x)}dx\geq \int_{\Omega} |u(x)|^{p_n^+}dx$$ and, thanks to \eqref{interm3}, we get that
 $$ 1= \lim_{n\to \infty} (\L^N(\Omega))^ { 1/p_n^+} \geq \lim_{n\to \infty} \bigg( \int_{\Omega} |u(x)| ^{p_n(x)}dx\bigg)^{ 1/p_n^+}\geq   \lim_{n\to \infty} \bigg( \int_{\Omega} |u(x)|^{p_n^+}dx\bigg)^{ 1/p_n^+}=||u||_{\infty}=1$$
 that is 
 \beq\label{interm}
\lim_{n\to \infty} \bigg( \int_{\Omega} |u(x)| ^{p_n(x)}dx\bigg)^{ 1/p_n^+}=1.
\eeq

 Due to \eqref{pn2},  the sequence $\beta_n=\left(\frac{p_n^+}{p_n^-} \right)_n$ satisfies  $1\leq \beta_n\leq \beta$. Then \eqref{interm} implies
  \beq\label{interm1}\lim_{n\to \infty}  \bigg( \int_{\Omega} |u(x)| ^{p_n(x)}dx\bigg)^{ 1/p_n^- }=\lim_{n\to \infty} \bigg( \int_{\Omega} |u(x)| ^{p_n(x)}dx\bigg)^{  \beta_n/p_n^+}=1.\eeq
 Moreover, by \eqref{relaztotale},
we have  that 
 for every $n\in \N$ 
$$\min\bigg\{ \bigg( \int_{\Omega} |u(x)| ^{p_n(x)}dx\bigg) ^{1/p_n^- }, \bigg( \int_{\Omega} |u(x)| ^{p_n(x)}dx\bigg) ^{1/p_n^+ } \bigg\}\leq  ||u||_{p_n(\cdot)} $$
$$\leq \max\bigg\{ \bigg( \int_{\Omega} |u(x)| ^{p_n(x)}dx\bigg) ^{1/p_n^- }, \bigg( \int_{\Omega} |u(x)| ^{p_n(x)}dx\bigg) ^{1/p_n^+ }\bigg\}.
$$
Then, taking into account \eqref{interm} and  \eqref{interm1}, when we pass to the limit when $n\to\infty$,  we get
  $$ \lim_{n\to \infty}  ||u||_{p_n(\cdot)}=1.$$
\\
\\
\fbox{$(ii) \Longrightarrow (i) $}   Assume now that $\lim_{n\to \infty} ||u||_{p_n(\cdot)}\in \R$. Let $q\geq 1$. Thanks to \eqref{pn1}, there exists   $n_0=n_0(q)\in \N$ big enough such that $p_n^->q$ for every $n\geq n_0$.

By Lemma \ref{pmenoq}, we get that 

$$ ||u||_{q}\leq \max\bigg\{\big(\L^N(\Omega) \big)^{\frac {1} {q} -\frac {1} {p_n^-} } , \big(\L^N(\Omega)\big)^{\beta \big (\frac {1} {q}-\frac {1} {p_n^+}\big)}\bigg\}  \bigg[ 1+\frac q {p_n^{+}(\cdot)} (\beta- 1)\bigg ]^{1/q}||u||_ {p_n(\cdot)}$$ 
for every $n\geq n_0$.

Since for every $q\geq 1$ we have that  \beq\label{limit}\lim_{n\to \infty}\bigg[ 1+\frac q {p_n^{+}(\cdot)} (\beta- 1)\bigg ]^{1/q}=1 \eeq
by passing to the limit when $n\to \infty$ it follows that 
 \begin{eqnarray*}  ||u||_{q}&\leq&  \lim_{n\to \infty} \max\bigg\{\big(\L^N(\Omega)\big)^{\frac {1} {q} -\frac {1} {p_n^-} } , \big(\L^N(\Omega)\big)^{\beta \big (\frac {1} {q}-\frac {1} {p_n^+}\big) } \bigg \}  \bigg[ 1+\frac q {p_n^{+}(\cdot)} (\beta- 1)\bigg ]^{1/q}||u||_ {p_n(\cdot)}\\
&=& \max\big\{\big(\L^N(\Omega) \big)^{\frac {1} {q}  } , \big(\L^N(\Omega)\big)^{ \frac {\beta} {q} }  \} \lim_{n\to \infty}  ||u||_ {p_n(\cdot)}\in \R \quad \forall q\geq 1 . \end{eqnarray*}
This implies \begin{eqnarray*}\lim_{q\to \infty}  ||u||_{q}\leq
 \lim_{q\to \infty} \left [
\max\big\{\big(\L^N(\Omega) \big)^{\frac {1} {q}  } , \big(\L^N(\Omega)\big)^{ \frac {\beta} {q} }  \} \lim_{n\to \infty} ||u||_ {p_n(\cdot)} \right ]\le 
 \lim_{n\to \infty}||u||_ {p_n(\cdot)}\in \R.\end{eqnarray*}
Then $u\in L^\infty(\Omega)$ and by the first part of this proof, it holds $||u||_{\infty}=\lim_{n\to \infty} ||u||_{p_n(\cdot)}.$ 

\qed

We conclude this section with  the following lemma, already shown in \cite{AP} when  $f=f(x,\xi)$ is a Carath\'eodory integrand (see Lemma 4.5 therein). For the reader's convenience we report here the proof.

\begin{lemma}\label{liminfFp} Let $f:\Om\times\R^d\times \R^{k} \to\R^+$ be a normal integrand. Then   
$$
\lim_ {q\to \infty} \left( \int_\Omega \int_{\R^k}
f(x,v(x),\xi)^q d\mu_x (\xi) dx\right)^{1/q} = \supess_{x \in \om} \left(
\mu_x\hbox{-}\supess_{\xi\in\R^k} f(x,v(x),\xi) \right),
$$
for every  Young measure $\mu=(\mu_x)$ and for every measurable function $v:\Omega\to \R^d$. 

\end{lemma}
{\bf Proof.} Taking into account Theorem \ref{FTYM}, part (1),  the following inequality 
\begin{eqnarray*}
& & \limsup_ {q\to \infty} \left( \int_\Omega \int_{\R^{k}}
f(x,v(x),\xi)^q d\mu_x (\xi) dx\right)^{1/q}\\
&\leq & \limsup_ {q\to \infty} \left( \int_\Omega  \mu_x(\R^{k}) \mu_x\hbox{-}\supess_{\xi\in\R^{k}} (f(x,v(x),\xi))^q dx\right)^{1/q}\\
&\leq & \limsup_ {q\to \infty} \left( \int_\Omega \left(\mu_x\hbox{-}\supess_{\xi\in\R^{k}} f(x,v(x),\xi) \right)^q dx\right)^{1/q}\\
&\leq &  \supess_{x \in \om} \left(
\mu_x\hbox{-}\supess_{\xi\in \R^{k}} f(x,v(x),\xi) \right)
\end{eqnarray*}
is straighforward, by the convergence of the $L^q$ norms to the $L^{\infty}$ norm. Let us prove that
$$
\liminf_ {q\to \infty} \left( \int_\Omega \int_{\R^{k}}
f(x,v(x),\xi)^q d\mu_x (\xi) dx\right)^{1/q} \geq  \supess_{x \in \om} \left(
\mu_x\hbox{-}\supess_{\xi\in\R^{k}} f(x,v(x),\xi) \right).
$$
 Without loss of generality we assume that
\beq\label{ineq}
\liminf_ {q\to \infty} \left( \int_\Omega \int_{\R^{k}}
f(x,v(x),\xi)^q d\mu_x (\xi) dx\right)^{1/q}<+\infty.
\eeq
For every fixed exponent $r$ such that $q>r$, by applying H\"older's  inequality   we get that
\begin{equation}\label{qrJ}
\left( \int_\Omega \int_{\R^{k}} f(x,v(x),\xi)^q d\mu_x (\xi)
dx\right)^{1/q} \geq  \left( \int_\Omega \left( \int_{\R^{k}} f(x,v(x),\xi)^r d\mu_x (\xi) \right)^{q/r} dx\right)^{1/q}.
\end{equation}
Passing to the limit as $q\to\infty$, by the convergence of the
$L^q$-norm to the $L^\infty$-norm, we have that
\begin{equation}\label{qress}
\lim_{q\to\infty}\left( \int_\Omega \left( \int_{\R^{k}}
f(x,v(x),\xi)^r d\mu_x (\xi) \right)^{q/r} dx\right)^{1/q}= \supess_{x \in
\om} \left( \int_{\R^{k}} f(x,v(x),\xi)^r d\mu_x (\xi)
\right)^{1/r}\,.
\end{equation}
We now denote 
$$
g_r(x):=\left( \int_{\R^{k}} f(x,v(x),\xi)^r d\mu_x (\xi)
\right)^{1/r}.
$$ 
Then $(g_r)$ is an increasing positive family  pointwise  converging  to the function 
$$g(x):=\mu_x\hbox{-}\supess_{\xi\in\R^{k} } f(x,v(x),\xi)$$ as $r\to\infty$.
Moreover, by (\ref{ineq})-(\ref{qress}), we have that $\sup_r ||g_r||_{\infty}<+\infty$.
In particular, by Lebesgue's dominated convergence theorem, we have that $g_r\wto g$ weakly* in $L^{\infty}$. By   (\ref{qrJ}), (\ref{qress}) and the weak* lower semicontinuity of the $L^\infty$-norm, we have that
$$
\liminf_ {q\to \infty} \left( \int_\Omega \int_{\R^{k}}
f(x,v(x),\xi)^q d\mu_x (\xi) dx\right)^{1/q} \ge
\supess_{x \in \om}\left(\mu_x\hbox{-}\supess_{\xi\in\R^{k}} f(x,v(x),\xi) \right),
$$
which concludes the proof. \qed

\section{The $L^p$ approximation }\label{Lp}

In this section we study the $L^p$-approximation, via $\Gamma$-convergence, of supremal functionals.
In the following we consider  a sequence $(p_n)$ of functions $p_n: \overline{\Omega} \rightarrow (1, + \infty) $, satisfying \eqref{pn1} and \eqref{pn2}  and a normal integrand  $f:\Om\times\R^d\times\R^{Nd}\to\R$ satisfying the following assumptions:
\begin{enumerate} 
\item [{\bf (H1)}] for a.e. $x\in \Omega,$  $f(x,u,\cdot)$ is level convex for every $u\in \R^d$;
\item  [{\bf (H2)}] there exist $\alpha, \gamma>0$ such that
\beq\label{crescitabasso} f(x,u, \xi) \geq  \alpha |\xi|^{\gamma}  
\qquad \hbox{ for a.e } x\in \Om, \hbox{ for every }
(u, \xi)\in\R^d\times  \R^{Nd}.
\eeq
\end{enumerate} 

\subsection{Statement of the main results}


We start by stating all theorems to easily compare the  results obtained according to the different set of  hypotheses and topologies considered.

The following result requires a regularity assumption of $\Omega$ in the proof of the $\Gamma$-liminf inequality since  we use the Sobolev imbedding, but we drop the hypothesis  that $\Omega$ is connected (used in the proof  given in \cite{BM} when the authors use the Poincarè-Wirtinger inequality).

\begin{thm}\label{gamma} Let $\Omega\subseteq \R^N$ be a bounded open set with Lipschitz boundary.
Let $f:\Om\times\R^d\times\R^{Nd}\to\R$ be a normal integrand  satisfying  assumptions {\bf (H1)} and {\bf (H2)}. 
 Let $F_n:L^{1}(\Omega,\R^{d})\to [0, + \infty]$ be the functional defined by 
\beq\label{approxseq}
F_n(u):=\left\{\begin {array}{cl} \displaystyle \ || f(\cdot,u(\cdot),Du(\cdot))||_{p_n(\cdot)}
&  \hbox{if } \, u\in W^{1,p_n(\cdot)}(\Omega,\R^{d})\\
+\infty  & \hbox{otherwise,}
\end{array}\right.
\eeq
and let $F:L^{1}(\Omega,\R^{d})\to [0, + \infty]$ be the functional defined by 
\beq\label{myfunc} F(u):=\left\{\begin {array}{cl} \displaystyle
\supess_{\Omega} f(x,u(x),Du(x))
&  \hbox{if } \, u\in W^{1,\infty}(\Omega,\R^{d}),\\
+\infty  & \hbox{otherwise.}
\end{array}\right.
\end{equation}
    
Then, 
 \begin{enumerate}
 \item[(i)] for every $u\in L^{1}(\Omega,\R^{d})$ and $(u_n)\subset  L^{1}(\Omega,\R^{d})$ such that $u_n\wto u$  in $L^{1}(\Omega,\R^{d})$, we have
 $$F(u)\leq \liminf_{n\to\infty} F_n(u_n) ;$$
 \item [(ii)]  for every $u\in L^{1}(\Omega,\R^{d})$ there exists  $(u_n)\subset  L^{1}(\Omega,\R^{d})$
 such that $u_n\to u$  in $L^{1}(\Omega,\R^{d})$ and
 $$\limsup_{n\to\infty} F_n(u_n)\leq F(u).$$
 \end{enumerate}
 In particular,  $(F_n)$  $\Gamma\hbox{-}$ converges to $F$,  as $n\to +\infty$,  with
 respect to the $L^{1}$-{\bf strong  convergence}.
 \end{thm}
The following result instead  does not require any regularity assumption of $\Omega.$ 
 \begin{thm}\label{gamma2} Let $\Omega\subseteq \R^N$ be a bounded open set.
Let $f:\Om\times\R^d\times\R^{Nd}\to\R$ be a normal integrand satisfying assumptions {\bf (H1)} and {\bf (H2)}. 
Let $X\in \{ L^{\infty}(\Omega,\R^{d}), C(\Omega,\R^{d})\}$ be endowed with the norm $||\cdot||_{\infty}$. Let $F_n:X\to [0, + \infty]$ be the functional defined by 
\beq\label{approxseq}
F_n(u):=\left\{\begin {array}{cl} \displaystyle \ || f(\cdot,u(\cdot),Du(\cdot))||_{p_n(\cdot)}
&  \hbox{if } \, u\in W^{1,p_n(\cdot)}(\Omega, \mathbb{R}^d )\\
+\infty  & \hbox{otherwise,}
\end{array}\right.
\eeq
and let $F:X\to [0, + \infty]$ be the functional defined by 
\beq\label{myfunc} F(u):=\left\{\begin {array}{cl} \displaystyle
\supess_{\Omega} f(x,u(x),Du(x))
&  \hbox{if } \, u\in W^{1,\infty}(\Omega,\R^{d}),\\
+\infty  & \hbox{otherwise.}
\end{array}\right.
\end{equation}
    
Then, 
 \begin{enumerate}
 \item[(i)] for every $u\in X$ and $(u_n)\subset  X$ such that
 $u_n\to u$   in $X$, we have
 $$F(u)\leq \liminf_{n\to\infty} F_n(u_n) ;$$
 \item [(ii)]  for every $u\in X$ there exists  $(u_n)\subset  X$
 such that $u_n\to u$  in $X$ and
 $$\limsup_{n\to\infty} F_n(u_n)\leq F(u).$$
 \end{enumerate}
 In particular,  $(F_n)$  $\Gamma\hbox{-}$ converges to $F$,  as $n\to +\infty$,  with
 respect to the $L^{\infty}$-{\bf strong  convergence}.
 \end{thm}

As a corollary, by applying the previous result when   $(p_n)_n$ is an arbitrary real sequence diverging to $+\infty$,   we get the following improvement of Theorem 3.1 in \cite{CDP}.

\begin{cor}\label{gamma4} Let $\Omega\subseteq \R^N$ be a bounded open set.
Let $f:\Om\times\R^d\times\R^{Nd}\to\R$ be a normal integrand  satisfying assumptions {\bf (H1)} and {\bf (H2)}. 
 Let $X\in \{ L^{\infty}(\Omega,\R^{d}), C(\Omega,\R^{d})\}$ be endowed with the norm $||\cdot||_{\infty}$. For every $p\geq 1$ let  $F_p:X\to [0, + \infty]$ be the functional defined by \beq\label{approxseq1}
F_p(u):=\left\{\begin {array}{cl} \displaystyle \ || f(\cdot,u(\cdot),Du(\cdot))||_{p}
&  \hbox{if } \, u\in W^{1,p}(\Omega,\R^{d})\\
+\infty  & \hbox{otherwise,}
\end{array}\right.
\eeq
and let $F$ be the functional defined by \eqref{myfunc}. Then,   $(F_p)$  $\Gamma\hbox{-}$converges to $F$,  as $p\to +\infty$,  with
 respect to the uniform  convergence.
 \end{cor}
Finally we show the following results:
\begin{thm}\label{gamma5} Let $\Omega\subseteq \R^N$ be a bounded open set  with Lipschitz boundary.
Let $f:\Om\times\R^d\times\R^{Nd}\to\R$ be a normal integrand satisfying  assumptions {\bf (H1)} and {\bf (H2)}. 
Let $\F_n:L^{1}(\Omega,\R^{d})\to  [0, + \infty]$ be the functional defined by 
\beq\label{approxseqfin}
\F_n(u):=\left\{\begin {array}{cl} \displaystyle \ \int_{\Omega} \frac{1}{p_n(x)}  f^{p_n(x)}(x,u(x),Du(x))dx
&  \hbox{if } \, u\in W^{1,p_n(\cdot)}(\Omega,\R^{d})\\
+\infty  & \hbox{otherwise,}
\end{array}\right.
\eeq
and let $\F:L^{1}(\Omega,\R^{d})\to [0, + \infty]$ be the functional defined by 
\beq\label{myfuncfin} \F(u):=\left\{\begin {array}{cl} \displaystyle 0
&  \hbox{if } \, u\in W^{1,\infty}(\Omega,\R^{d}) \hbox{ and }||f(x,u(x),Du(x))||_{\infty}\leq 1,\\
+\infty  & \hbox{otherwise.}
\end{array}\right.
\end{equation}

Then, 
 \begin{enumerate}
 \item[(i)] for every $u\in L^{1}(\Omega,\R^{d})$ and $(u_n)\subset  L^{1}(\Omega,\R^{d})$ such that $u_n\wto u$  in $L^{1}(\Omega,\R^{d})$, we have
 $$\F(u)\leq \liminf_{n\to\infty} \F_n(u_n) ;$$
 \item [(ii)]  for every $u\in L^{1}(\Omega,\R^{d})$ there exists  $(u_n)\subset  L^{1}(\Omega,\R^{d})$
 such that $u_n\to u$  in $L^{1}(\Omega,\R^{d})$ and
 $$\limsup_{n\to\infty} \F_n(u_n)\leq \F(u).$$
 \end{enumerate}
 In particular,  $(\F_n)$  $\Gamma\hbox{-}$ converges to $\F$,  as $n\to +\infty$,  with
 respect to the $L^{1}$-{\bf strong  convergence}.
 \end{thm}

\begin{thm}\label{gamma6} Let $\Omega\subseteq \R^N$ be a bounded open set.
Let $f:\Om\times\R^d\times\R^{Nd}\to\R$ be a normal integrand satisfying  assumptions {\bf (H1)} and {\bf (H2)}. 
Let $X\in \{ L^{\infty}(\Omega,\R^{d}), C(\Omega,\R^{d})\}$ be endowed with the norm $||\cdot||_{\infty}$. Let $\F_n:X\to [0, + \infty]$ be the functional defined by \eqref{approxseqfin}
and let $\F:X\to [0, + \infty]$ be the functional defined by \eqref{myfuncfin}.
Then, 
 \begin{enumerate}
 \item[(i)] for every $u\in X$ and $(u_n)\subset  X$ such that
 $u_n\to u$   in $X$, we have
 $$\F(u)\leq \liminf_{n\to\infty} \F_n(u_n) ;$$
 \item [(ii)]  for every $u\in X$ there exists  $(u_n)\subset  X$
 such that $u_n\to u$  in $X$ and
 $$\limsup_{n\to\infty} \F_n(u_n)\leq \F(u).$$
 \end{enumerate}
 In particular,  $(\F_n)$  $\Gamma\hbox{-}$ converges to $\F$,  as $n\to +\infty$,  with
 respect to the $L^{\infty}$-{\bf strong  convergence}.
 \end{thm}

\subsection{Proofs of {Theorems}}
%


\noindent {\bf Proof of Theorem \ref{gamma}.}  First of all we consider  the case when $\gamma\geq 1$.  We observe that
\begin{equation}\label{gamlimsup}
\limsup_ {n\to \infty}\ F_n(u)  \leq F(u)
\end{equation}
for any $u \in L^{1}(\om,\R^{d})$. Indeed, if $F(u)=+\infty$, there
is nothing to prove, and if $F(u)< +\infty$, then,  $u\in W^{1,\infty}(\om,\R^d) $ and $f(\cdot,u(\cdot), Du(\cdot))\in \li(\Omega)$. By Proposition \ref{infty} we have that
$$\lim_ {n\to \infty}F_n(u)
=||f(\cdot,u(\cdot), Du(\cdot))||_{\infty}=F(u) $$ so that
(\ref{gamlimsup}) follows. As a consequence, for any $u \in L^{1}(\om,\R^{d})$ it holds
$$ 
\Gamma(L^1) \hbox{-}\limsup_{n\to \infty} F_n (u)\leq \limsup_{n\to \infty} F_n (u)
\leq F(u)\,.
$$
We now deal with the  liminf inequality. Let $(u_n)\in L^{1}(\om,\R^d)$ be a
sequence weakly converging to $u\in
L^{1}(\om,\R^d)$. 
Without loss of generality, we can assume that \beq\label{C}\liminf_{n\to
\infty}F_n(u_n) =\lim_{n\to
\infty}F_n(u_n)=M<+\infty,\eeq
hence, by definition of the functionals $F_n$, we have that there exists $n_0\in \N$ such that 
$u_n\in W^{1,p_n(\cdot)}(\om,\R^d)$ for every $n\geq n_0$. 
Fix $q> 1$ and let  $n_1\geq n_0$ be  such that, in view of \eqref{pn1} 
\begin{equation}
\label{pallino}
p_n^-\geq q \hbox{ and } F_n(u_n)\leq M+1 \qquad \forall\, n\geq n_1.
\end{equation}
Then, by applying Lemma \ref{pmenoq} to $u_n$ and to $Du_n$ with $p(\cdot)=p_n(\cdot)$, we get that $u_n\in W^{1,q}(\om,\R^{d})$  for every $n\geq n_1;$  on the other hand, still be Lemma \ref{pmenoq}, we also get the estimate 
\begin{eqnarray} \label{lp} 
 ||u_n||_{q}\leq   c_{q,n}||u_n||_ {p_n(\cdot)}  \qquad \forall n \geq n_1,
 \end{eqnarray}
where $$c_{q,n}:=\max  \left\{\big (\L^N(\Omega) \big)^{\frac {1} {q} -\frac {1} {p_n^-} } , \big(\L^N(\Omega)\big)^{\beta \big (\frac {1} {q}-\frac {1} {p_n^+}\big) }  \right\}  \bigg[ 1+\frac q {p_n^{+}(\cdot)} (\beta- 1)\bigg ]^{1/q}.
$$
Moreover  the function $v_n(\cdot):=f(x,u_n(\cdot),Du_n(\cdot)) \in L^{p_n(\cdot)}(\om)$  $\forall\, n\geq n_1 $ and, by applying again Lemma \ref{pmenoq},  this time  to $v_n$, we obtain 
that for every $n\geq n_1$
\begin{eqnarray}\label{un}   ||f(x,u_n(\cdot),Du_n(\cdot)) ||_{q}
\label{fun} &\leq& c_{q,n} ||f(x,u_n(\cdot),Du_n(\cdot))||_ {p_n(\cdot)}\\
&\leq &  (M+1)c_{q,n}, \nonumber
\end{eqnarray}
 where we used \eqref{pallino}.
Note that  for every fixed $q>1$ the sequence $(c_{q,n})_n$ is bounded  since,  by \eqref{pn1}  $$\lim_{n\to \infty}c_{q,n}= \max  \left \{ \big(\L^N(\Omega) \big)^{\frac {1} {q}} , \big(\L^N(\Omega)\big)^{\frac {\beta } {q}}   \right \}:=c_q<+\infty.$$
Taking into account the growth condition (\ref{crescitabasso}),   \eqref{un}  implies that for every $n\geq n_1$ 
%
\begin{eqnarray*} \label{grad} ||Du_n||^{\gamma}_{ q} & \leq& \big(\L^N(\Omega) \big)^{\frac {\gamma} q-\frac 1  { q}} ||Du_n||^{\gamma}_{\gamma q} \\ 
&\leq& \big(\L^N(\Omega) \big)^{\frac {\gamma} q-\frac 1  { q}}   \frac {1} {\alpha}  ||f(x,u_n(\cdot),Du_n(\cdot)) ||_{q}  \leq \big(\L^N(\Omega) \big)^{\frac {\gamma} q-\frac 1  { q}}  \frac {M+1} {\alpha}     c_{q,n}.
 \end{eqnarray*}
 that is 
 \beq \label{gradlimq}  ||Du_n||_{ q} \leq  \big(\L^N(\Omega) \big)^{\frac 1  q-\frac 1  { {\gamma} q} }  \left(\frac {M+1} {\alpha}  c_{q,n}\right)^{\frac 1 \gamma}.
 \eeq
In particular 
\begin{eqnarray} \label{gradlim}  \sup_{n\geq n_1}  ||Du_n||_{q}\leq  \big(\L^N(\Omega) \big)^{\frac 1  q-\frac 1  { {\gamma} q}}  \left(\frac {M+1} {\alpha}  \sup_{n\geq n_1}  c_{q,n}\right)^{\frac 1 \gamma} <+\infty.\end{eqnarray}

Then, up to a subsequence (depending on $q$), $(Du_n)_n$ weakly converges to a function $w$ in $L^q(\Omega,\R^{Nd})$. Since $(u_n)_n$  weakly converges to $u$ in $L^1(\Omega,\R^d)$,  it is easy to show  that $w$ is the distributional gradient of $u$. In particular  $u\in W^{1,1}(\om,\R^d)$ and, since every subsequence  of $(Du_n)_n$ admits a subsequence converging to $Du$, we get that  the whole sequence  $Du_n\wto Du$ weakly in $L^q(\Omega,\R^{Nd})$.
  Now we show that $u\in W^{1,q}(\om,\R^d)$.
  Note that, being  $u\in W^{1,1} (\om,\R^d)$, thanks to the Sobolev immersion, we get that $u\in L^{1^*}(\om,\R^d)=L^{\frac N {N-1}}(\om,\R^d).$
 Since $Du\in L^{\frac N {N-1}}(\Omega,\R^{Nd})$, we deduce that $u\in W^{1,\frac N {N-1}}(\om,\R^d).$
 Then $u\in L^{(\frac N {N-1})^*}(\om,\R^d)=L^{\frac N {N-2}}(\om,\R^d).$
 By going on, after $k=N-1$ steps we get that 
 $$u\in L^{(\frac N {N-(k-1)})^*}(\om,\R^d)= L^{\frac N {N-k}}(\om,\R^d)= L^{N}(\om,\R^d)$$ that is $u\in W^{1,N}(\om,\R^d)$. By Sobolev immersion, we can conclude that $u\in L^q(\om,\R^d)$ for every $q\geq N$ and, since $Du\in L^q(\om,\R^{Nd})$ for every $q\geq 1$,  we obtain that $u\in W^{1,q}(\om,\R^d)$ for every $q\geq 1$
and $u_n\wto u$ weakly in $W^{1,q}(\om,\R^d)$.  In particular $u\in L^{\infty}(\Omega,\R^d)$.
%

%
%
 %

Moreover, taking into account \eqref{grad}, we get     
\begin{eqnarray*} ||Du||_{q}\leq \liminf_{n\to \infty}  ||Du_n||_{q}\leq \big(\L^N(\Omega) \big)^{\frac 1  q-\frac 1  { {\gamma} q} }  \left( \frac {M+1} {\alpha}    c_q \right)^{\frac 1 \gamma}    \  \ \forall q>1\end{eqnarray*}  
that implies, taking into account that $c_{q}\to 1$ when $q\to \infty$, \begin{equation}\label{stimaunif}  \lim_{q\to \infty} ||Du||_{q}\leq  \left( \frac {M+1} {\alpha} \right)^{\frac 1 \gamma} <\infty \end{equation}
i.e.  $Du\in L^{\infty}(\om,\R^{Nd})$ and $u\in W^{1,\infty}(\om,\R^d)$.

%
%
%

By Remarks  \ref{limitedebole} and   \ref{rem2.3}, we have that $(Du_n)$ generates a Young measure
$(\mu _x)_{x\in\Omega}$ such that $\mu_x(\R^{Nd})=1$ and
\beq\label{gradiente}
Du(x)=\int_{\R^{Nd}}\xi \,d\mu _{x} (\xi)
\eeq
for a.e. $x\in\Omega$.
 Then,  for any fixed $q>N$, by applying (\ref{un}) and \eqref {limsotto},   we have
that
\begin{eqnarray*} \liminf_{n\to \infty} F_n(u_n)&\geq &   \liminf_{n \rightarrow + \infty}  \frac 1 {c_{q,n}} ||f(\cdot,u_n(\cdot), Du_n(\cdot)||_{q} \\
 &=& \frac 1 {c_{q}} \liminf_{n\to \infty}
\left( \int_\Omega f^q(x,u_n(x),Du_n(x)) dx \right)^{1/q} \\
& \geq & \left( \int_\Omega \int_{ \R^{Nd}} f^q(x,u(x), \xi) d\mu_x
(\xi) dx\right)^{1/q}.
\end{eqnarray*}
By applying Lemma \ref{liminfFp} we obtain 
\begin{eqnarray}\label{liminf1} \liminf_{n\to \infty} F_n(u_n)
 & \geq &  \liminf_{q\to \infty} \left( \int_\Omega \int_{\R^{Nd}} f^q(x,u(x),\xi) d\mu_x
(\xi) dx\right)^{1/q}\\
\nonumber  & = &\supess_{x \in \om} \left( \mu_x\hbox{-}\supess_{\xi\in\R^{Nd}} f(x,u(x),\xi)
\right).
\end{eqnarray}
Since $f(x,u(x),\cdot)$ is  level convex for a.e. $x\in\Omega$, taking into account \eqref{gradiente}, by Jensen's inequality \eqref {JI} we have that 
$$
f(x,u(x),Du(x))=f \left (x,u(x),\int_{\R^{Nd}}\xi \,d\mu _{x} (\xi)\right)\leq \mu_{x}\hbox{-}\supess_{\xi\in\R^{Nd}}  f(x,u(x), \xi )
$$
for a.e. $x\in\Omega$. In particular
\begin{equation}\label{sup1}
\supess_{x\in\Omega}f(x,u(x),Du(x))\leq  \supess_{x\in \Om}
\Bigl(\mu_{x}\hbox{-}\supess_{\xi\in\R^{Nd}}   f(x,u(x),\xi)\Bigr)\,.
\end{equation}
%
%
Then, by  the very definition of $F$,  we get 
$$
 F(u)=\supess_{x\in\Omega}f(x,u(x),Du(x))
$$
and    \eqref{sup1} and \eqref{liminf1} 
imply   the $\Gamma$-liminf inequality. 

 Thus, the proof in the case $\gamma\geq 1$ is concluded.
 Assume now that $0<\gamma<1$.  First of all we observe that, since the function $t\to t^{\frac 1 \gamma}$ is monotone on $[0,+\infty)$,  then  the function $g(x,u,\xi):=f^{\frac 1 \gamma}(x,u,\xi)$ is level convex too with respect to the gradient variable and satisfies the grouth condition 
$$ g(x,u, \xi) \geq  \alpha' |\xi|$$
 for a.e $ x\in \Om$,  for every $(u, \xi)\in\R^d\times  \R^{Nd}
$, with   $\alpha'=\alpha^{\frac 1 \gamma}$.

Then, we get that
the sequence of the functionals $G_n:L^{1}(\Omega,\R^{d})\to \R\cup \{+\infty\}$ defined by 
\beq\label{approxseqGn}
G_n(u):=\left\{\begin {array}{cl} \displaystyle \ || g (\cdot,u(\cdot),Du(\cdot))||_{\gamma p_n(\cdot)}
&  \hbox{if } \, u\in W^{1,\gamma p_n}(\Omega,\R^{d})\\
+\infty  & \hbox{otherwise,}
\end{array}\right.
\eeq
$\Gamma$-converges to 
 $G:L^{1}(\Omega,\R^{d})\to \R\cup \{+\infty\}$  defined by 
\beq\label{myfuncG} G(u):=\left\{\begin {array}{cl} \displaystyle
\supess_{\Omega} g(x,u(x),Du(x))
&  \hbox{if } \, u\in W^{1,\infty}(\Omega,\R^{d}),\\
+\infty  & \hbox{otherwise,}
\end{array}\right.
\end{equation}
with
 respect to the $L^{1}$- strong  convergence.
Since $p_n^-\to +\infty$, for $n$ big enough we have that $\frac 1 \gamma<p_n^-$. Then, by Lemma \ref{lemma3bocea} applied with $s=\frac 1 \gamma$, we get that
\beq\label{bocea2} || g (\cdot,u(\cdot),Du(\cdot))||_{\gamma p_n(\cdot)}^{ \gamma}=|| f^{\frac 1 \gamma} (\cdot,u(\cdot),Du(\cdot))||_{\gamma p_n(\cdot)}^{ \gamma}
= || f (\cdot,u(\cdot),Du(\cdot))||_{ p_n(\cdot)}.
\eeq
%
Moreover, since $\gamma<1$, we have that  $$W^{1, n}(\Omega,\R^{d})\subseteq W^{1, \gamma n}(\Omega,\R^{d}).$$
Thus, taking into account  \eqref{bocea2}, we get
$$G_n^{\gamma} \leq F_n\leq F.$$
By passing to the $\Gamma$-limit when $n\to \infty$ with respect to the $L^1$-convergence and noticing that $G_n^\gamma$ $\Gamma$-converges to 
$G^{\gamma}=F$, we get  the thesis.   \QED

%

\noindent {\bf Proof of Theorem \ref{gamma2}.}  As in the proof of Theorem \ref{gamma} it is sufficient to prove the result in the case $\gamma\geq 1$. The proof of the $\Gamma$-limsup inequality  follows the same arguments as in Theorem \ref{gamma}. In order to get the  $\Gamma$-liminf inequality,  it is sufficient to note that 
if $(u_n)\subseteq X$ is a
sequence  $L^\infty$-converging to $u$ in $X$, then  $(u_n)$  weakly $L^q$-converges to $u$ for every $q\geq 1$. By applying  inequality  \eqref{gradlim} we get that the sequence  $(Du_n)_n$ weakly converges to $Du$ in $L^{q}(\Omega,\R^d)$ for every $q>1$. 
 In particular $(u_n)_n$ converges weakly to $u$ in $W^{1,q}(\Omega,\R^{Nd})$ for every $q>N$. Then $(Du_n)$ generates a Young measure
$(\mu _x)_{x\in\Omega}$ such that $\mu_x(\R^{Nd})=1$  and 
$Du(x)=\int_{\R^{dN}}\xi \,d\mu _{x} (\xi)
$
for a.e. $x\in\Omega$. Thus the $\Gamma$-liminf inequality follows by applying Jensen's inequality \eqref {JI} and Lemma \ref{liminfFp}  in order to get   \eqref{sup1}.  \QED
 
\noindent {\bf Proof of Corollary \ref{gamma4}.} It is sufficient to apply Theorem \ref{gamma2} to get that,
for every  sequence $(p_n)_n$ diverging  to 
$\infty$ as $n\to\infty$, the sequence  $(F_n)$, defined by \eqref{approxseq},
 $\Gamma\hbox{-}$converges to $F$ with
 respect to the uniform  convergence.
 \QED
\noindent {\bf Proof of Theorem \ref{gamma5}.}  We observe that
\begin{equation}\label{gamlimsup1}
\limsup_ {n\to \infty}\ \F_n(u)  \leq \F(u)
\end{equation}
for any $u \in L^{1}(\om,\R^{d})$. Indeed, if $\F(u)=+\infty$, there
is nothing to prove, and if $\F(u)< +\infty$, then $\F(u)=0$  that implies  $u\in W^{1,\infty}(\om,\R^d) $ and $||f(\cdot,u(\cdot), Du(\cdot))||_{\infty}\leq 1$.
In particular $$
\limsup_{n\to \infty}\F_n(u)= \limsup_{n\to \infty}\int_{\Omega} \frac{1}{p_n(x)}  f^{p_n(x)}(\cdot,u(\cdot),Du(\cdot))dx \leq  \L^N(\Omega)   \limsup_{n\to \infty} \frac{1 }{p^-_n} = 0.$$
Then it is sufficient to take  $u_n=u$ to get the $\Gamma$-limsup inequality.
We now deal with the $\Gamma$-liminf inequality. Let $(u_n)\in L^{1}(\om,\R^d)$ be a
sequence weakly converging to $u\in
L^{1}(\om,\R^d)$. 
Without loss of generality, we can assume that \beq\label{C1}\liminf_{n\to
\infty}\F_n(u_n) =\lim_{n\to
\infty}\F_n(u_n)=M<+\infty,\eeq

hence, by definition of the functionals $\F_n$, we have that there exists $n_0\in \N$ such that $\F_n(u_n)\leq 2M$ for every $n\geq n_0$. In particular 
$u_n\in W^{1,p_n(\cdot)}(\om,\R^d)$ for every $n\geq n_0$. 
 For each $n\in \N$, define $$ \Omega^+_{n}:=\{ x\in \Omega: f(x,u_n(x),Du_n(x))> 1\} \qquad \textnormal{and} \qquad\Omega^-_{n}:=\{x\in \Omega: f(x,u_n(x),Du_n(x))\leq 1\}.$$
Then, for every $n\geq n_0$ it holds
\begin{eqnarray*}
&& \frac{1}{p_n^+}   \int_{\Omega}  f^{p_n(x)}(x,u_n(x),Du_n(x))dx \\
&\le& \int_{\Omega} \frac{1}{p_n(x)}  f^{p_n(x)}(x,u_n(x),Du_n(x))dx\\
&= &\int_{\Omega^+_n}   \frac{1}{p_n(x)}  f^{p_n(x)}(x,u_n(x),Du_n(x))dx +\int_{\Omega_n^-}  \frac{1}{p_n(x)}  f^{p_n(x)}(x,u_n(x),Du_n(x)) dx\\
&\leq &  \int_{\Omega^+_n}   \frac{1}{p_n(x)}  f^{p_n(x)}(x,u_n(x),Du_n(x))dx +\L^N(\Omega)  \frac{1}{p_n^-} \\
&\leq &  \int_{\Omega}   \frac{1}{p_n(x)}  f^{p_n(x)}(x,u_n(x),Du_n(x))dx +\L^N(\Omega)  \frac{1}{p_n^-} \\
&\leq & 2M +\L^N(\Omega)  \frac{1}{p_n^-}.
\end{eqnarray*}

In particular
$$  \int_{\Omega}  f^{p_n(x)}(x,u_n(x),Du_n(x))dx \leq p_n^+  \left (2M +\L^N(\Omega)  \frac{1}{p_n^-} \right)  $$
that implies
$$  \left [ \rho_{p_n(\cdot)}( f(x,u_n(\cdot),Du_n(\cdot)))\right ]^{\frac 1 {p_n^+}}   \leq  \left [p_n^+  \left(2M +\L^N(\Omega)  \frac{1}{p_n^-}\right)\right]^{\frac 1 {p_n^+} } = \left [2M p_n^+  +\L^N(\Omega)  \frac{p_n^+} {p_n^-} \right]^{\frac 1 {p_n^+} }=M(n)$$
and  also, replacing $p_n^+$ with $p_n^-$ in the exponent of the modular in the left hand side
$$  \left [ \rho_{p_n(\cdot)}( f(x,u_n(\cdot),Du_n(\cdot))) \right ]^{\frac 1 {p_n^-}}   \leq  \left [p_n^+  \left(2M +\L^N(\Omega)  \frac{1}{p_n^-}\right) \right ]^{\frac 1 {p_n^-} } =\left( M(n) \right)^{  \frac {p_n^+} {p_n^-}  } .$$
Taking into account \eqref{relaztotale},
it follows

 \begin{equation} \label{hn}  \|f(\cdot,u_n(\cdot),Du_n)\|_{p_n(\cdot)} \leq \max  \left \{  M(n),  \left( M(n) \right)^{  \frac {p_n^+} {p_n^-}}   \right \}
\end{equation}
Since $(\frac  {p_n^+}  {p_n^-})_n$ is a bounded sequence, $M(n)\to 1$ and $\frac 1 {p_n^+}\to 0$  when $n\to \infty$, the previous inequality implies that the  sequence $ (\|f(\cdot,u_n(\cdot),Du_n)\|_{p_n(\cdot)})_n $ is bounded and 
$$\liminf_{n\to \infty} \|f(\cdot,u_n(\cdot),Du_n)\|_{p_n(\cdot)}\leq 1.$$
Moreover, by applying the $\Gamma$-liminf inequality in Theorem \ref{gamma}, we obtain that $u\in \wi$ and 
$$F(u)\leq \liminf_{n\to \infty} F_n(u_n)=\liminf_{n\to \infty} \|f(\cdot,u_n(\cdot),Du_n)\|_{p_n(\cdot)}\leq 1$$
where $F$ is defined  by \eqref{myfunc}.
This implies $\F(u)=0$. \QED

\noindent {\bf Proof of Theorem \ref{gamma6}.} The proofs follows the lines of the previous result by applying Theorem \ref{gamma2} instead of  Theorem \ref{gamma}.\QED

%

\noindent{\bf Acknowledgments.}  
Both the authors are members of GNAMPA-INdAM, whose support is gratefully acknowledged.   The work of the ME is also supported by the University of Modena and Reggio Emilia through the project FAR 2019 ''Equazioni differenziali: problemi evolutivi, variazionali ed applicazioni'', Coord. Prof. Maria Manfredini. ME is indebted with Dipartimento di Matematica and Informatica of University of Ferrara for its kind support and hospitality.

\end{document}